\documentclass{article}

\usepackage{amscd}
\usepackage{amsmath}
\usepackage{amssymb}
\usepackage{color}
\usepackage{graphicx}
\usepackage{latexsym}
\usepackage{makeidx}
\usepackage{mathtools}
\usepackage{theorem}
\usepackage{xypic}

\usepackage[colorlinks]{hyperref}
\hypersetup{colorlinks=false , linkcolor=blue}

\newtheorem{theorem}{Theorem}[section]

\newtheorem{lemma}[theorem]{Lemma}
\newtheorem{corollary}[theorem]{Corollary}

\DeclareMathOperator{\Ker}{Ker}

\def\d{\mathrm{d}}
\def\raya{\ \underline{\phantom{a}}\ }

\def\qed{\hspace*{\fill }$\square $ }

\begin{document}

\author{Alberto Navarro \thanks{Instituto de Ciencias Matem\'{a}ticas (CSIC-UAM-UC3M-UCM)}}

\title{On Grothendieck's Riemann-Roch Theorem}

\date{March 12$^{\mathrm{th}}$, 2016}

\maketitle

\begin{abstract}
We prove that, for smooth quasi-projective varieties over a field, the $K$-theory
$K(X)$ of vector bundles is the universal cohomology theory where $c_1(L\otimes \bar
L)=c_1(L)+c_1(\bar L)-c_1(L)c_1(\bar L)$. Then, we show that Grothendieck's
Riemann-Roch theorem is a direct consequence of this universal property, as well as
the universal property of the graded $K$-theory $GK^\bullet (X)\otimes \mathbb{Q}$.
\end{abstract}
\bigskip

\noindent \textbf{INTRODUCTION:} In 1957 Grothendieck introduced the $K$-group $K(X)$
of vector bundles on an algebraic variety $X$ and the $K$-group $G(X)$ of coherent
sheaves. He showed that the $K$-theory $K(X)$ has a cohomological behavior whereas
$G(X)$ has a homological behavior. For example, $K(X)=G(X)$ when $X$ is smooth, so
that we have an unexpected direct image $f_!\colon K(Y)\to K(X)$ for any projective
map $f\colon Y\to X$ between smooth varieties. Then Grothendieck gave his astonishing
formulation of the Riemann-Roch theorem as a determination of the relationship
between this direct image $f_!$ on $K$-theory and the usual direct image $f_*$ at the
level of the Chow ring $CH^\bullet(X)$ (or the singular cohomology, etc.). The
relationship involved two series: the Chern character $\mathrm{ch}$ and the Todd
class $\mathrm{Td}$.

\bigskip

\noindent \textbf{Grothendieck's Riemann-Roch Theorem:} \emph{Let $f\colon Y\to X$ be
a projective morphism between smooth quasi-projective algebraic varieties and denote
$T_X$ and $T_Y$ their tangent bundles. Then we have a commutative square}
$$
\xymatrix{K(Y)  \ar[r]^-{f_!} \ar[d]_{\mathrm{Td}(T_Y)\cdot
\mathrm{ch}} & K(X) \ar[d]^{\mathrm{Td}(T_X)\cdot \mathrm{ch}} \\
CH^\bullet(Y)\otimes \mathbb{Q} \ar[r]^-{f_*} &CH^\bullet(X)\otimes \mathbb{Q}. }
$$

\bigskip

In this article we prove that Grothendieck's $K$-theory $K(X)$ is the universal
cohomology theory where Chern classes of line bundles fulfil that
$$
c_1(L\otimes \bar L)=c_1(L)+c_1(\bar L)-c_1(L)c_1(\bar L).
$$
Even if this universal property is clearly implicit in the work of Levine, Morel,
Panin, Pimenov and R\"{o}ndigs on the algebraic cobordism (\cite{Levine},\cite{Panin3}),
this result seems to be new.

After that, we show that Grothendieck's Riemann-Roch theorem is a direct consequence
of this universal property of the $K$-theory. We also obtain the universal property
of the graded $K$-theory $GK^\bullet (X)\otimes \mathbb{Q}$, when the ring $K(X)$ is
filtered by the subgroups $F^d(X)$ generated by the coherent sheaves with support of
codimension $\ge d$.

\medskip

Let us comment the plan of this article. In section 1 we review the theory of Chern
classes in a cohomology theory $A(X)$. In section 2 we prove the universal property
of $K$-theory: for any cohomology theory $A(X)$ following the multiplicative law
$c_1(L\otimes \bar L)=c_1(L)+c_1(\bar L)-c_1(L)c_1(\bar L)$ there exists a unique
functorial ring morphism $\varphi \colon K(X)\to A(X)$ preserving direct images. The
key point, due to Panin in \cite{Panin}, is to prove that any functorial ring
morphism preserving Chern classes of the tautological line bundles also preserves
direct images,
$\varphi (f_!(y))=f_*(\varphi (y))$. For the sake of completeness, we include Panin's
proof, greatly simplified in the case of the projection $p\colon \mathbb{P}^n\to
\mathrm{pt}$ using the approach of D\'{e}glise (\cite{Deglise2}). For $p$, we observe
that the inverse image $A(\mathrm{pt})\to A(\mathbb{P}^n)$ is dual to the direct
image $A(\mathbb{P}^n)\to A(\mathrm{pt})$ with respect to the bilinear pairing
defined by the fundamental class of the diagonal
$$
\Delta _n\in A(\mathbb{P}^n\times \mathbb{P}^n)=A(\mathbb{P}^n)\otimes
_{A(\mathrm{pt})}A(\mathbb{P}^n).
$$
Since $\varphi $ preserves inverse images by hypothesis and it also preserves the
fundamental class of the diagonal (because it preserves the direct image by the
diagonal immersion), the proof is reduced to check that the metric of the diagonal is
non-singular.

In addition, we compute all possible direct images on a given cohomology theory.
We show that, given a formal series $F(t)=1+a_1t+a_2t^2+\ldots $ with coefficients
in the
cohomology ring $A(\mathrm{pt})$ of a point, we may modify the direct image so that
the new Chern class of a line bundle $L\to X$ is $c_1^{\mathrm{new}}(L)=c_1(L)\cdot
F(c_1(L))$.

In section 3 we consider cohomology theories $A(X)$ following the additive law
$c_1(L\otimes \bar L)=c_1(L)+c_1(\bar L)$. In any such a theory we may modify the
direct image with a formal series with rational coefficients (essentially an
exponential) so that it follows the multiplicative law of the $K$-theory and hence,
due to the universal property of $K$-theory, we have a functorial ring morphism
$\mathrm{ch}\colon K(X)\to A(X)\otimes \mathbb{Q}$ compatible with the new direct
image. This compatibility is just Grothendieck's Riemann-Roch theorem. We also show
that the graded $K$-theory $GK^\bullet (X)\otimes \mathbb{Q}$ is the universal graded
cohomology with rational coefficients (with values in the category of graded
$\mathbb{Q}$-algebras). Finally, in section 4 we present some direct consequences of
these theorems.

\medskip

\noindent \textbf{Acknowledgments:} I would like to thank my advisors, Fr\'{e}d\'{e}ric
D\'{e}glise and Jos\'{e} Ignacio Burgos Gil, for generously and patiently guiding me and
sharing with me their knowledge during the preparation of this work.
I am deeply thankful to my father Juan Antonio, to whom I dedicate this paper, among
many things for our mathematical discussions about the Riemann-Roch theorem.

The author has been partially supported by ICMAT Severo Ochoa project SEV-2015-0554
(MINECO) and MTM2013-42135-P (MINECO).

\section{Cohomology theories}

Throughout this paper, let us fix a base field $k$.

\medskip

\noindent \textbf{Definition:}  With Panin (\cite{Panin}), we define a
\textbf{cohomology theory} to be a contravariant functor $A$ from the category of
smooth quasi-projective varieties over $k$ into the category of commutative rings,
endowed with a functorial morphism of $A(X)$-modules $f_*\colon A(Y)\to A(X)$, called
\textbf{direct image}, for any projective morphism $f\colon Y\to X$ (that is to say,
$\mathrm{Id}_*=\mathrm{Id}$, $(fg)_*=f_*g_*$ and the projection formula
$f_*(f^*(x)y)=xf_*(y)$ holds). Hence we have a \textbf{fundamental class}
$[Y]^A\coloneqq i_*(1)\in A(X)$ for any  smooth closed subvariety $i\colon Y\to X$,
and a \textbf{Chern class} $c_1^A(L)\coloneqq s_0^*(s_{0*}(1))\in A(X)$ for any line
bundle $L\to X$ (where $s_0\colon X\to L$ is the zero section).

These data are assumed to satisfy the following conditions:
\begin{enumerate}
\item
The ring morphism $i_1^*+i_2^*\colon A(X_1\amalg X_2)\to A(X_1)\oplus A(X_2)$ is an
isomorphism, where $i_j\colon X_j\to X_1\amalg X_2$ is the natural immersion.

Hence $A(\emptyset )=0$.

\item
The ring morphism $\pi ^*\colon A(X)\to A(P)$ is an isomorphism for any affine bundle
$\pi \colon P\to X$ (i.e., a torsor over a vector bundle $E\to X$).

Therefore, if $\pi \colon L\to X$ is a line bundle, then $c_1^A(L)=s^*(s_{0*}(1))$
for any section $s$, because $s^*=s_0^*$, both being the inverse of $\pi ^*\colon
A(X)\to A(L)$.

\item
For any smooth closed subvariety $i\colon Y\to X$ we have an exact sequence $\
A(Y)\xrightarrow{i_*}A(X)\xrightarrow{j^*}A(X-Y)\ $ .

\item
If a morphism $f\colon \bar X\to X$ is transversal to a smooth closed subvariety
$i\colon Y\to X$ of codimension $d$ (that is to say, $f^{-1}(Y)=\emptyset$ or $\bar
Y=Y\times _X\bar X$ is a smooth subvariety of $\bar X$ of codimension $d$, so that
the natural epimorphism $f^*N^*_{Y/X}\to N^*_{\bar Y/\bar X}$ is an isomorphism,
where $N_{Y/X}$ denotes the normal bundle), then we have a commutative square
$$
\xymatrix{A(Y) \ar[r]^-{f^*} \ar[d]^{i _*} & A(\bar Y) \ar[d]^{i _*}\\ A(X)
\ar[r]^-{f^*} & A(\bar X). }
$$

\item
Let $\pi \colon \mathbb{P}(E)\to X$ be a projective bundle. For any morphism $f\colon
Y\to X$ we have a commutative square
$$
\xymatrix{A(\mathbb{P}(E))\ar[r]^-{f^*} \ar[d]^{\pi _*} & A(\mathbb{P}(f^*E))
\ar[d]^{\pi _*}\\ A(X) \ar[r]^-{f^*} & A(Y). }
$$

\item
Let $\pi \colon \mathbb{P}(E)\to X$ be the projective bundle associated to a vector
bundle $E\to X$ of rank $r+1$ and let $\xi _E\to \mathbb{P}(E)$ be the tautological
line bundle. Consider the structure of $A(X)$-module in $A(\mathbb{P}(E))$ defined by
the ring morphism $\pi ^*\colon A(X)\to A(\mathbb{P}(E))$, and put $x_E=c_1^A(\xi
_E)$, then $1,x_E,\ldots ,x_E^r$ define a basis :
$$
A(\mathbb{P}(E))=A(X)\oplus A(X)x_E\oplus \ldots \oplus A(X)x_E^r.
$$
\end{enumerate}

Given cohomology theories $A,\, \bar A$ on the smooth quasi-projective $k$-varieties,
a \textbf{morphism} of cohomology theories $\varphi \colon A\to \bar A$ is a natural
transformation preserving direct images. That is to say, for any smooth
quasi-projective variety $X$ we have a ring morphism $\varphi \colon A(X)\to \bar
A(X)$ such that $\varphi (f^*(a))=f^*(\varphi (a))$, $a\in A(X)$, for any morphism
$f\colon Y\to X$, and $\varphi (f_*(b))=f_*(\varphi (b))$, $b\in A(Y)$, for any
projective morphism $f\colon Y\to X$.\bigskip

\goodbreak

\noindent \textbf{Remarks:}

\begin{enumerate}

\item
Both $i_1$ and $i_2$ are transversal to $i_1$ and $i_2$; hence $i_2^*i_{1*}=0$,
$i_1^*i_{2*}=0$, $i_1^*i_{1*}=\mathrm{Id}$, $i_2^*i_{2*}=\mathrm{Id}$, so that the
inverse of the isomorphism $i_1^*+i_2^*$ of axiom 1 is just $i_{1*}+i_{2*}\colon
A(X_1)\oplus A(X_2)\to A(X_1\amalg X_2)$, and for any projective morphism
$f=f_1\amalg f_2\colon X_1\amalg X_2\to Z$ we have that $f_*=f_{1*}+f_{2*}\colon
A(X_1)\oplus A(X_2)\to A(Z)$.

\item
Let $E\to X$ be a vector bundle of rank $r$. In $\mathbb{P}(E)$ we have an exact
sequence $0\to \xi _E\to \pi ^*E\to Q\to 0$, where $Q$ is a vector bundle of rank
$r-1$, and the ring morphism $\pi ^*\colon A(X)\to A(\mathbb{P}(E))$ is injective by
axiom 6. Proceeding by induction on the rank we see that the following
\textbf{splitting principle} holds: \emph{There exists a base change $\pi \colon
X'\to X$ such that $\pi ^*E$ admits a filtration $0=E_0\subset E_1\subset \ldots
\subset E_{r}=\pi ^* E$ whose quotients $E_i/E_{i-1}$, $i=1,\ldots ,r$, are line
bundles, and $\pi ^*\colon A(X)\to A(X')$ is injective}.

\item
Given a line bundle $L\to X$ and a morphism $f\colon Y\to X$, the induced morphism
$f\times 1 \colon f^*L=Y\times _XL\to X\times _XL=L$ is transversal to the null
section $s_0\colon X\to L$. Hence
$$
f^*s_0^*s_{0*}(1)=\bar s_0^*(f\times 1)^*s_{0*}(1)=\bar s_0^*\bar s_{0*}(1),
$$
where $\bar s_0$ stands for the null section of $f^*L$, and we see that Chern classes
are functorial,
\begin{equation}\label{1}
c_1(f^*L)=f^*c_1(L).
\end{equation}

\item
Let $L_Y\to X$ be the line bundle defined by a smooth closed hypersurface $i\colon
Y\to X$ (dual to the line bundle defined by the ideal of $Y$). It admits a section
$s\colon X\to L_Y$ vanishing just on $Y$ and transversal to the null section $s_0$.
Hence
\begin{equation*}\label{2}
c_1(L_Y)=s^*(s_{0*}(1))=i_*(i^*(1))=i_*(1)=[Y]\in A(X).
\end{equation*}
In particular, if $\xi _d\to \mathbb{P}^d$ is the tautological line bundle of the
projective space of dimension $d$, the Chern class of the dual bundle $\xi _d^*$ is
just the fundamental class of an hyperplane, $c_1(\xi _d^*)=[\mathbb{P}^{d-1}]$.

In general, $L_x$ will denote a line bundle with Chern class $x\in A(X)$, and we say
that a cohomology theory $A$ follows the additive group law $x+y$ when
$c_1^A(L_x\otimes L_y)= x + y$, that it follows the multiplicative\footnote{Since
$(1-x)(1-y)=1-(x+y-xy)$, if we declare $x$ to be the coordinate of $1-x$ in the
multiplicative group $k-\{ 0\} $, then the group law is just $x+y-xy$.} group law
$x+y-xy$ when $c_1^A(L_x\otimes L_y)= x+y-xy$, and so on.

\end{enumerate}

\goodbreak

\noindent \textbf{Examples:}

\begin{enumerate}

\item
The $K$-theory $K(X)$ is a cohomology theory (\cite{SGA6}, \cite{Weibel}). The
fundamental class of a smooth closed subvariety $i\colon Y\to X$ is
$$
[Y]^K=i_!(1)=\mathcal{O}_Y\in K(X),
$$
where $\mathcal{O}_Y$ is the structural sheaf of $Y$. The Chern class of a line
bundle $L\to X$ is
$$
c_1^K(L)=s^!_0(s_{0!}(1))=1-L^*\in K(X).
$$

In general $c_1^K(L^*)\neq -c_1^K(L)$; but we have
$$
1-(L\otimes \bar L)^*=(1-L^*)+(1-\bar L^*)-(1-L^*)(1-\bar L^*),
$$
so that the $K$-theory follows the multiplicative group law $x+y-xy$,
$$
c_1^K(L\otimes \bar L)=c_1^K(L)+c_1^K(\bar L)-c_1^K(L)c_1^K(\bar L).
$$

\item
Denote $F^d(X)$ the subgroup of $K(X)$ generated by the coherent sheaves with support
of codimension $\ge d$. The graded $K$-theory $GK^\bullet (X)=\oplus
_dF^d(X)/F^{d+1}(X)$ is a cohomology theory (\cite{SGA6}, \cite{Navarro}).
The fundamental class of a smooth closed subvariety $i\colon Y\to X$ of codimension
$d$ is
$$
[Y]^{GK}=[\mathcal{O}_Y]\in F^d(X)/F^{d+1}(X),
$$
and the Chern class of a line sheaf $L$ is
$$
c_1^{GK}(L)=[1-L^*]=[L-1]\in F^1(X)/F^2(X).
$$

Since $(1-L^*)(1-\bar L^*)\in F^2(X)$, the graded $K$-theory follows the additive law
$c_1^{GK}(L\otimes \bar L)=c_1^{GK}(L)+c_1^{GK}(\bar L)$.

\item
When $k=\mathbb{C}$, the singular cohomology classes of even degree define a
cohomology theory $H^{2\bullet }(X,\mathbb{Z})=\oplus _dH^{2d}(X,\mathbb{Z})$
following the additive law $x+y$.

\item
The Chow ring $CH^\bullet (X)$ of rational equivalence classes of cycles is a
cohomology theory (\cite{Fulton}) also following the additive law $x+y$. The
fundamental class of a subvariety Z is given by the rational equivalence class $[Z]$.

\item
Given a field extension $k\to K$ and a cohomology theory $A$ on the smooth
quasi-projective $K$-varieties, then $X\mapsto A(X\times _kK)$ is a cohomology theory
on the smooth quasi-projective $k$-varieties. Hence, if $k$ is a field with a fixed
embedding $k\to \mathbb{C}$, then we have a cohomology theory $H^{2\bullet }(X\times
_k\mathbb{C},\mathbb{Z})$ on the smooth quasi-projective $k$-varieties.

\end{enumerate}

\noindent \textbf{Definition:} Let $E\to X$ be a vector bundle of rank $r$. With
Grothendieck (\cite{Groth}) we define the \textbf{Chern classes} $c_n^A(E)\in A(X)$
of $E$ to be the coefficients of the characteristic polynomial
$c(E)=x^r-c_1^A(E)x^{r-1}+\ldots +(-1)^rc_r^A(E)$ of the endomorphism of the free
$A(X)$-module $A(\mathbb{P}(E))$ defined by the multiplication by $x_E=c_1^A(\xi
_E)$,
\begin{equation}\label{3}
x_E^r-c_1^A(E)x_E^{r-1}+\ldots +(-1)^rc_r^A(E)=0,
\end{equation}
and we put $c_n(E)$ when the cohomology theory $A$ is clear.

\medskip

\noindent \textbf{Remark:} The signs are introduced so that the Chern class of a line
bundle $L$ coincides with the former one. In fact, $\mathbb{P}(L)=X$ and $\xi _L=L$,
so that $x_L=c_1(\xi _L)=c_1(L)$.

\begin{theorem}\label{4}
For any morphism $f\colon Y\to X$, we have $c_n(f^*E)=f^*\bigl( c_n(E)\bigr) $.
\end{theorem}

\noindent \emph{Proof:} Any morphism $f\colon Y\to X$ induces a morphism $f\colon
\mathbb{P}(f^*E)\to \mathbb{P}(E)$ such that $f^*\xi _E=\xi _{f^*E}$. Hence
$f^*(x_E)=x_{f^*E}$ by (\ref{1}). Applying $f^*$ to (\ref{3}), we have
$$
x_{f^*E}^r-\bigl( f^*c_1(E)\bigr) x_{f^*E}^{r-1}+\ldots +(-1)^rf^*c_r(E)=0,
$$
and we conclude that $c_n(f^*E)=f^*c_n(E)$. \qed

\begin{theorem}
Chern classes are additive. That is to say, $c(E)=c(E_1)\cdot c(E_2)$ for any exact
sequence $0\to E_1\to E\to E_2\to 0$ of vector bundles,
$$
c_n(E)=\sum _{i+j=n}c_i(E_1)\cdot c_j(E_2)\qquad ;\quad n,i,j\in \mathbb{N}.
$$
\end{theorem}

\noindent \emph{Proof:} Assume that $E_1$ is a line bundle. Then $i\colon
X=\mathbb{P}(E_1)\to \mathbb{P}(E)$ is a section of $\mathbb{P}(E)\to X$, so that
$i_*$ is injective. Moreover $i^*\xi _E=\xi _{E_1}$, $i^*(x_E)=x_{E_1}$.

Denote $j\colon U\to \mathbb{P}(E)$ the complement of $\mathbb{P}(E_1)$. The natural
projection $p\colon U\to \mathbb{P}(E_2)$ is an affine bundle (of associated vector
bundle $\mathrm{Hom}(\xi _{E_2},p^*E_1)$) and $j^*\xi _E=p^*\xi _{E_2}$, so that
$j^*(x_E^n)=p^*(x_{E_2}^n)$.

Hence $j^*\colon A(\mathbb{P}(E))\to A(U)\stackrel{p^*}{=}A(\mathbb{P}(E_2))$ is
surjective and, by axiom 3, we have a commutative diagram with exact rows (the first
square commutes by the projection formula)
$$
\xymatrix{0\ar[r] & A(\mathbb{P}(E_1)) \ar[d]^{\cdot  x_{E_1}} \ar[r]^{i_*} &
A(\mathbb{P}(E))\ar[r]^{j^*}
 \ar[d]^{\cdot x_E} & A(\mathbb{P}(E_2)) \ar[r] \ar[d]^{\cdot  x_{E_2}} & 0\\
0\ar[r] & A(\mathbb{P}(E_1))\ar[r]^{i_*} & A(\mathbb{P}(E))\ar[r]^{j^*} &
A(\mathbb{P}(E_2)) \ar[r] & 0}
$$

In the general case, by the splitting principle, we may assume that we have a line
bundle $L\subset E_1$ such that $\bar E_1=E_1/L$ and $\bar E=E/L$ are vector bundles,
so that we have an exact sequence $0\to \bar E_1\to \bar E\to E_2\to 0$ and we
conclude by induction on the rank of $E_1$,
$$
c(E)=c(L)c(\bar E)=c(L)c(\bar E_1)c(E_2)=c(E_1)c(E_2).
$$
\qed

\medskip

\noindent \textbf{Remark:} By the above theorem $1+c_1^A(E)t+c_2^A(E)t^2+\ldots
+c_r^A(E)t^r$ is an additive function on the vector bundles over $X$, with values in
the multiplicative group of invertible formal series with coefficients in $A(X)$.
Hence it extends to the $K$-group, and we obtain Chern classes $c_n^A\colon K(X)\to
A(X)$.

Now, given a vector bundle $E\to X$ of rank $r$, by the splitting principle there is
a base change $\pi \colon X'\to X$ such that $\pi ^*\colon A(X)\to A(X')$ is
injective and $\pi ^*E=L_{\alpha _1}+\ldots +L_{\alpha _r}$ is a sum of line bundles
in $K(X')$. Therefore, $c_n(E)$ is the $n$-th elementary symmetric function of the
"roots" $\alpha _1,\ldots ,\alpha _r$,
\begin{equation}\label{roots}
c_n(E)=\sum _{i_1<\ldots <i_n}\alpha _{i_1}\ldots \alpha _{i_n}.
\end{equation}

For example, in the $K$-theory the Chern class of a line bundle $L$ is just
$c_1^K(L)=1-L^*$; hence the first Chern class of a vector bundle $E$ of rank $r$ is
\begin{equation}\label{5}
c_1^K(E)=(1-L^*_{\alpha _1})+\ldots +(1-L^*_{\alpha _r})=r-E^*,
\end{equation}
and $c_r^K(E)=(1-L^*_{\alpha _1})\cdot \ldots \cdot (1-L^*_{\alpha _r})=\sum
_i(-1)^i\Lambda ^iE^*$.

\begin{corollary}\label{Coro fibrado proyectivo}
The cohomology ring of the projective spaces $\mathbb{P}^d$ is
$$
A(\mathbb{P}^d)=A(\mathrm{pt})[x]/(x^{d+1})=A(\mathrm{pt})[y]/(y^{d+1}),
$$
where $x$ corresponds to $x_d=c_1(\xi _d)$ and $y$ corresponds to $y_d=c_1(\xi
_d^*)=[\mathbb{P}^{d-1}]$.
\end{corollary}

\noindent \emph{Proof:} The projective space $\mathbb{P}^d$ is just the projective
bundle of a trivial vector bundle of rank $d+1$ over a point. By additivity, trivial
bundles have null Chern classes; hence $x^{d+1}_d=0$ in
$A(\mathbb{P}^d)=A(\mathrm{pt})\oplus A(\mathrm{pt})x_d\oplus \ldots \oplus
A(\mathrm{pt})x_d^d$.

Now, in $\mathbb{P}^1$ we have an exact sequence $0\to \xi _1\to 1\oplus 1\to \xi
_1^*\to 0$, where $1\oplus 1$ stands for the trivial vector bundle of rank 2. Hence
$y_1=-x_1$ in $A(\mathbb{P}^1)$. Considering a line $\mathbb{P}^1\to \mathbb{P}^d$,
we see that $y_d=-x_d+a_2x^2_d+\ldots +a_dx^d_d$ in $A(\mathbb{P}^d)$. Hence
$y^{d+1}_d=0$ and we conclude. \qed

\begin{corollary}\label{Coro Chern nilpotente}
Chern classes are always nilpotent.
\end{corollary}

\noindent \emph{Proof:} Let $L\to X$ be a line bundle. By Jouanolou's trick
(\cite{Jouanolou}) there is an affine bundle $p\colon P\to X$ such that $P$ is an
affine variety.

Now $p^*L$ is generated by global sections; hence $p^*L=f^*(\xi ^*_d)$ for some
morphism $f\colon P\to \mathbb{P}^d$. It follows that $p^*c_1(L)=f^*(y_d)$ is
nilpotent, since so is $y_d$, and we see that $c_1(L)$ is nilpotent, because $p^*$ is
an isomorphism. We conclude since any Chern class, after a base change injective in
cohomology, is a sum of products of Chern classes of line bundles.\qed

\bigskip

\section{Universal Property of the $K$-theory}

\begin{theorem}
If a cohomology theory $A$ follows the group law $x+y-xy$ of the $K$-theory, there is
a unique morphism of cohomology theories $\varphi \colon K\to A$.
\end{theorem}

\noindent \emph{Proof:} Let $E\to X$ be a vector bundle. Due to (\ref{5}), we have
$E=\mathrm{rk}\, E-c_1^K(E^*)$ in $K(X)$; hence the unique possible morphism $\varphi
\colon K\to A$ is
\begin{equation}\label{10}
\varphi (E)\coloneqq \mathrm{rk}\, E-c_1^A(E^*).
\end{equation}

Now, since $\mathrm{rk}$ and $c_1^A$ are additive, $\varphi $ is an additive function
on the vector bundles over $X$. Therefore $\varphi $ defines a group morphism
$\varphi \colon K(X)\to A(X)$, and $\varphi $ commutes with inverse images because so
do the rank, and $c_1^A$ by (\ref{4}).

This group morphism $\varphi $ preserves products of line bundles because $A$ follows
the law $x+y-xy$,
\begin{align*}
\varphi (L_1\otimes  L_2)&=1-c_1^A(L_1^*\otimes L_2)=1-c_1^A(L_1^*)-c_1^A(L_2^*)+
c_1^A(L_1^*)c_1^A( L_2^*) \\
&=\bigl( 1-c_1^A(L_1^*)\bigr) \bigl( 1-c_1^A(L_2^*)\bigr) =\varphi (L_1)\cdot \varphi
(L_2).
\end{align*}
Hence $\varphi $ is a ring morphism, $\varphi (ab)=\varphi (a)\varphi (b)$, since we
may assume (by the splitting principle) that $a,b\in K(X)$ are sums and differences
of line bundles.

It is only left for us to prove that $\varphi $ preserves direct images. It preserves
Chern classes of line bundles,
$$
\varphi \bigl( c_1^K(L)\bigr) =\varphi (1-L^*)=1-\varphi (L^*)=1-\bigl(
1-c_1^A(L)\bigr) =c_1^A(L).
$$
Therefore $\varphi $ preserves fundamental classes of hypersurfaces, and the theorem
follows from the following result. \qed

\bigskip

\noindent \textbf{Panin's Lemma:} (\cite{Panin}) \emph{Let $A$ and $\bar A$ be two
cohomology theories on smooth quasi-projective varieties. If a natural transformation
$\varphi \colon A\to \bar A$ preserves the first Chern class of  the tautological
line bundles $\xi _d\to \mathbb{P}^d$ (i.e. $\varphi (c_1^A(\xi_d))=c_1^{\bar
A}(\xi_d)$) then it preserves direct images:}
\begin{equation}\label{11}
\varphi (f_*(a))=\bar f_*(\varphi (a))
\end{equation}
\emph{for any projective morphism $f\colon Y\to X$, and any element $a\in
A(Y)$.}\bigskip

\goodbreak

\noindent \emph{Proof:} Let $Y$ be a hypersurface of a smooth quasi-projective
variety $X$. By the same argument as in Corollary \ref{Coro Chern nilpotente} we
deduce from Jouanolou's trick that the natural transformation $\varphi \colon A\to
\bar A$ preserves fundamental classes of hypersurfaces, $\varphi ([Y]^A)=[Y]^{\bar
A}$.

Now let $f\colon Y\to X$ be a projective morphism. By definition, $f$ is the
composition of a closed immersion $i\colon Y\to \mathbb{P}^n\times X$ with the
natural projection $\pi _X\colon \mathbb{P}^n\times X\to X$. If (\ref{11}) holds for
$i$ and $\pi _X$, then it also holds for the composition $f=\pi _X\circ i$, and it is
enough to prove (\ref{11}) for a closed immersion $i\colon Y\to X$ and the canonical
projection $\pi _X\colon \mathbb{P}^n\times X\to X$.
\begin{enumerate}
\item
\emph{If equation (\ref{11}) holds for the zero section $s\colon Y\to \bar
N=\mathbb{P}(1\oplus N_{Y/X})$ of the projective closure of the normal bundle
$N_{Y/X}$, then it also  holds for the closed immersion} $i\colon Y\to X$.

\emph{Proof:} Let $X'$ be the blow-up of $X\times \mathbb{A}^1$ along $Y\times 0$, so
that we have a commutative diagram
$$
\xymatrix{\bar N\ar[r]^-{i_0}& X'&X\ar[l]_-{i_1}\\
Y\ar[r]^-{i_0} \ar[u]^{s_0} &Y\times \mathbb{A}^1\ar[u]_{i} &Y\ar[l]_-{i_1}\ar[u]_{i}}
$$
Put $U=X'-(Y\times \mathbb{A}_1)$. By axiom 4 we have a commutative diagram
$$
\xymatrix{& \bar A(U)\\
\bar A(\bar N) & \bar A(X') \ar[l]_-{\bar i_0^*} \ar[u]_-{\bar j^*}\\
\bar A(Y) \ar[u]^-{\bar s_{0*}} & \bar A(Y\times \mathbb{A}_1) \ar[u]_-{\bar i_*}
\ar[l]_-{\bar i_0^*}^-{\sim } }
$$
Note that $(\mathrm{Ker}\, \bar i_0^*)\cap (\mathrm{Ker}\, \bar j^*)=0$. Indeed the
column is exact by axiom 3, and $\bar s_{0*}$ is injective because $s_0$ is a section
of the natural projection $\bar N\to Y$. Now consider the commutative diagram
$$
\xymatrix{&\bar A(U)&\\
\bar A(\bar N)&\bar A(X')\ar[l]_-{\bar i_0^*} \ar[u]_{\bar j^*} \ar[r]^-{\bar
i_1^*}&\bar A(X)\\
A(Y)\ar[u]^-{\Psi_1}&A(Y\times \mathbb{A}_1) \ar[l]_-{\sim } \ar[u]_-{\Psi _2}
\ar[r]_-{\sim } &A(Y) \ar[u]_-{\Psi _3}}
$$
where the vertical arrows are the difference of the morphisms that (\ref{11}) asserts
the coincidence ($\Psi _1=\bar s_{0*}\varphi -\varphi s_{0*}$ and so on). The
morphism $\Psi _1$ is zero by hypothesis, the morphism $\Psi _2$ is zero because
$\bar i_0^*\Psi _2=0$, $\bar j^*\Psi _2=0$ and $(\mathrm{Ker}\, \bar i_0^*)\cap
(\mathrm{Ker}\, \bar j^*)=0$, and we conclude that $0=\Psi _3=\bar i_*\varphi
-\varphi i_*$.

\item
\emph{Equation (\ref{11}) holds for the zero section $s\colon Y\to \bar
E=\mathbb{P}(1\oplus E)$ of the projective closure of any vector bundle $E\to Y$.}

\emph{Proof:} When $E=L$ is a line bundle, note that $\varphi (s_*(1))=\bar s_*(1)$
since $Y$ is an hypersurface in $\bar L$ and that $s^*\colon A(\bar L)\to A(Y)$ is
surjective. If  $a=s^*b\in A(Y)$, then equation (\ref{11}) holds:
$$
\varphi (s_*a)=\varphi (s_*s^*b)=\varphi (bs_*(1))=\varphi (b)\bar s_*(1)=\bar
s_*(\bar s^*\varphi (b))= \bar s_*(\varphi (a)).
$$
Therefore by the previous point it also holds for the immersion of any closed
hypersurface.

Now, if $E$ admits a filtration $\{ E_i\} $ such that the quotients $E_i/E_{i-1}$ are
line bundles, then (\ref{11}) holds for the zero section $Y\to \bar E_1$ and for the
morphisms $\bar E_1\to \bar E_2\to \ldots \to \bar E_r=\bar E$; hence it holds for
the composition $s\colon Y\to \bar E$.

In general, we have a morphism $\pi \colon Y'\to Y$ such that $\bar \pi ^*\colon \bar
A(Y)\to \bar A(Y')$ is injective and $E'=\pi ^*E$ admits such filtration. Then
(\ref{11}) holds for the zero section $s'\colon Y'\to \bar E'$, and we conclude
applying axiom 4 to the morphisms $\pi \colon \bar E'\to \bar E$ and $s\colon Y\to
\bar E$,
$$
\bar \pi ^*\bar s_*\varphi =\bar s'_*\bar \pi ^*\varphi =\bar s'_*\varphi \pi
^*=\varphi s'_*\pi ^*=\varphi \pi ^*s_*=\bar \pi ^*\varphi s_*.
$$

\item
\emph{If equation (\ref{11}) holds for the projection $p\colon \mathbb{P}^n\to
\mathrm{pt}$ onto a point, then it also holds does for the canonical projection} $\pi
_X\colon \mathbb{P}^n\times X\to X$.

\emph{Proof:} It follows from axiom 5.

\item
\emph{Equation (\ref{11}) holds for the projection $p\colon \mathbb{P}^n\to
\mathrm{pt}$ onto a point.}

Consider the closed immersion $i\colon \mathbb{P}^{n-1}\to \mathbb{P}^n$ and set
\begin{align*}
A&=A(\mathrm{pt}),\ x_n=c_1^A(\xi _n^*)=i_*(1)\in A(\mathbb{P}^n),\\
\bar A&=\bar A(\mathrm{pt}),\ \bar x_n=c_1^{\bar A}(\xi _n^*)=\bar i_*(1)\in \bar
A(\mathbb{P}^n).
\end{align*}
\indent By hypothesis $\varphi (x_n)=\bar x_n$; hence $\varphi (x_n^r)=\bar x_n^r$
and by axiom 6 the ring morphism $\varphi \colon A(\mathbb{P}^n)\to \bar
A(\mathbb{P}^n)$ induces an isomorphism of $\bar A$-algebras $A(\mathbb{P}^n)\otimes
_A\bar A=\bar A(\mathbb{P}^n)$.

We have to check that the $\bar A$-linear map $\bar p_*\colon \bar A(\mathbb{P}^n)\to
\bar A$ is obtained by base change of the $A$-linear map $p_*\colon
A(\mathbb{P}^n)\to A$.

Let $\Delta _n=\Delta _*(1)\in A(\mathbb{P}^n\times
\mathbb{P}^n)=A(\mathbb{P}^n)\otimes _AA(\mathbb{P}^n)$ be the fundamental class of
the diagonal immersion $\Delta \colon \mathbb{P}^n\to \mathbb{P}^n\times
\mathbb{P}^n$. We have
$$
(p_*\otimes 1)(\Delta _n)=(p\times 1)_*\Delta _*(1)=\mathrm{Id}_*(1)= 1\in
A(\mathbb{P}^n),
$$
where $p\times 1\colon \mathbb{P}^n\times \mathbb{P}^n\to \mathbb{P}^n$ is the second
projection. That is to say, $p_*$ corresponds to the unity, by means of the polarity
$A(\mathbb{P}^n)^*\to A(\mathbb{P}^n)$, $\omega \mapsto (\omega \otimes 1)(\Delta
_n)$, defined by the diagonal.

According to the next lemma, the linear map $p_*\colon A(\mathbb{P}^n)\to A$ is fully
determined by this condition. Since the fundamental class of the diagonal is stable
by the base change $A\to \bar A$, because (\ref{11}) holds for closed immersions, we
conclude that $p_*$ also is stable.\hfill $\square $

\end{enumerate}

\begin{lemma}
The metric $\Delta _n\in A(\mathbb{P}^n)\otimes _AA(\mathbb{P}^n)$ of the diagonal is
non-singular.
\end{lemma}

\noindent \emph{Proof:} By induction on $n$ we prove that
$$
\Delta _n=\sum _{r,s=0}^na_{rs}x_n^r\otimes x_n^s=\left( \raisebox{0.5\depth}{
\xymatrix@=1ex{ 0\ar@{.}[rr]\ar@{.}[dd]
 & & 0\ar@{.}[lldd]
 & 1\ar@{-}[lllddd] \\
 & & & \bullet\ar@{.}[lldd]\ar@{.}[dd] \\
0 & & & \\
1 & \bullet\ar@{.}[rr] & & \bullet } } \right)
$$
where $a_{rs}=0$ when $r+s<n$, and $a_{rs}=1$ when $r+s=n$. Indeed,
$$
i_*(x_{n-1}^r)=i_*i^*(x_n^r)=x_n^r\cdot i_*(1)=x_n^{r+1},
$$
and by axiom 4 we have that the fundamental class $(1\otimes i_*)(\Delta _{n-1})$ of
the diagonal of $\mathbb{P}^{n-1}$ in $\mathbb{P}^{n-1}\times \mathbb{P}^n$ is just
$(i^*\otimes 1)(\Delta _n)$. Note that
\begin{align*}
(i^*\otimes 1)(\Delta _n)&=\sum _{r=0}^{n-1}\sum _{s=0}^na_{rs}x_{n-1}^r\otimes
x_n^s,\\
(1\otimes i_*)(\Delta _{n-1})&=\sum _{r=0}^{n-1}\sum _{s=0}^{n-1}a'_{rs}x_{n-1}^r
\otimes x_n^{s+1},
\end{align*}
where $\Delta _{n-1}=\sum _{rs}a'_{rs}x_{n-1}^r\otimes x_{n-1}^s$.

By induction on $n$, we obtain the result for the coefficients $a_{rs}$, $r<n$.

By symmetry we also obtain it for $a_{rs}$, $s<n$, and we conclude.

\qed

\medskip

With Panin's lemma we can also compute the possible direct images on a given
cohomology theory.

Let us consider the "roots" $\alpha _1,\ldots ,\alpha _r$ of a vector bundle $E\to X$
of rank $r$ (\ref{roots}). For any formal power series $F(t)=\sum _na_nt^n\in
A(\mathrm{pt})[[t]]$ we put
$$
F_+(E)=F(\alpha _1)+\ldots +F(\alpha _r)\in A(X),
$$
where we consider $F(\alpha _1)+\ldots +F(\alpha _r)$ as a power series in the
elementary symmetric functions, which are just the Chern classes of $E$. It is well
defined since Chern classes are nilpotent, and it is an additive function on the
vector bundles over $X$, so defining a functorial group morphism $F_+\colon K(X)\to
A(X)$, named \textbf{additive extension} of $F$.

Analogously, let $F(t)=a_0+a_1t+\ldots $ be an invertible series. We put
$$
F_\times (E)=F(\alpha _1)\cdot \ldots \cdot F(\alpha _r)\in A(X)^*,
$$
and we obtain the \textbf{multiplicative extension} $F_\times \colon K(X)\to A(X)^*$
of $F$.

\begin{theorem}
Let $A$ be a cohomology theory, and $f_*^{\mathrm{new}}$ be another class of direct
images for $A$. For any projective morphism $f\colon Y\to X$ denote $T_f\coloneqq
T_Y-f^*T_X\in K(Y)$ the virtual relative tangent bundle. Then there exist an
invertible series $F(t)\in A(\mathrm{pt})[[t]]$ such that
$$
f_*^{\mathrm{new}}(a)=f_*(F_\times(T_f)^{-1}\cdot a).
$$
Moreover, every invertible series defines new direct images by the preceding
formula.
\end{theorem}
\emph{Proof:} From Corollary \ref{Coro fibrado proyectivo} we have that
$A(\mathbb{P}^d)=A(\mathbb{P}^d)=A(\mathrm{pt})[c_1(\xi_d)]/(c_1(\xi_d)^{d+1})$ for
all $d$. Therefore it follows that there exists a series $b_0+b_1t+\cdots\in
A(\mathrm{pt})[[t]]$ such that
$$
c_1^{\mathrm{new}}(\xi_d)=b_0+b_1\cdot c_1(\xi_d))+\cdots
$$
for all $d$. For $d=0$ we have $\mathbb{P}^d=\mathrm{pt}$ so that $b_0=0$. For $d=1$
consider the closed immersion $i\colon \mathrm{pt} \to \mathbb{P}^1$ with open
complement $j\colon \mathbb{A}^1\to \mathbb{P}^1$. There are exact sequences
$$
\xymatrix{ A(\mathrm{pt})\ar@<-.5ex>[r]_{i_*}\ar@<.5ex>[r]^{i_*^{\mathrm{new}}}&
A(\mathbb{P}^1)\ar[r]^{j^*}&A(\mathbb{A}^1)}
$$
where $\Ker j^*=A(\mathrm{pt})\cdot c_1(\xi_1)=A(\mathrm{pt})\cdot
c_1^{\mathrm{new}}(\xi_1)$. Since $c_1^{\mathrm{new}}(\xi_1)=b_1\cdot
c_1(\xi_1)$ we conclude that $b_1$ is invertible. We set
$$
F(t)=b_1+b_2t+\cdots.
$$
Note that $c_1^\mathrm{new}(L)=c_1(L)F(c_1(L))$.

For any projective morphism
$f\colon Y\to X$ consider the map $f_*(F_\times(T_f)^{-1}\raya ): A(Y)\to A(X)$.
The pair $(A,f_*(F_\times(T_f)^{-1}\raya ))$ is also a cohomology theory. Indeed, all
the axioms are easy to check, except the last one. Put $x=x_E\in A(\mathbb{P}(E))$
and $y=x_E^\mathrm{new}=xF(x)=x+\ldots \in A(\mathbb{P}(E))$, so that $y^n=x^n+\ldots
$, and we have $0=x^d=y^d$ for some exponent $d$. Since the powers of $x$ generate
the $A(X)$-module $A(\mathbb{P}(E))$, so do the
powers of $y$. Now, since $A(\mathbb{P}(E))$ is a free $A(X)$-module of rank $r+1$,
we conclude that $1,\ldots ,y^r$ define a basis (just consider the characteristic
polynomial of the endomorphism defined by $y$).

We can explicitly compute in $(A,f_*(F_\times(T_f)^{-1}\raya ))$
fundamental classes of a smooth closed hypersurface $i\colon Y\to X$:
$$
i_*\bigl( F_\times (N_{Y/X})\cdot 1\bigr) =i_*\bigl( F_\times
(i^*L_Y)\bigr) =F_\times (L_Y)i_*(1)=[Y]\cdot F([Y]).
$$
Therefore the first Chern class of line sheaf $L$ is
$s_0^*s_{0*}(F_\times (s_0^*L) \cdot 1)=c_1(L)F(c_1(L))$, because
$c_1(L)=s_0^*s_{0*}(1)$.

Consider the identity $A \to A$. Applying Panin's lemma to
$(A,f_*^\mathrm{new})$ and $(A,f_*(F_\times(T_f)^{-1}\raya ))$ we conclude.

\qed

\section{Riemann-Roch Theorem}

Let $A$ be cohomology theory following the additive law. When we consider formal
series with rational coefficients, the additive and multiplicative formal groups are
isomorphic, because of the exponential series, and we may modify the direct image of
$A\otimes \mathbb{Q}$ with an exponential so that the new theory follows the
multiplicative law $x+y-xy$ of the $K$-theory.

Since $e^{at}=1-(1-e^{at})$, we must fix a formal series $F(t)$ such that
\begin{equation}\label{8}
c_1^\mathrm{new}(L_x)=1-e^{ax}=xF(x).
\end{equation}

Now, $1-e^{at}=-at+\ldots $, so that it is convenient to fix $a=-1$.

Hence, so as to transform the additive law of $A \otimes \mathbb{Q}$ into a
multiplicative law, just modify the direct image with the formal series
\begin{equation}\label{9}
F(t)=\frac{1-e^{-t}}{t}=1-\frac{t}{2!}+\frac{t^2}{3!}+\ldots
\end{equation}
so that the new cohomology theory $A^{\mathrm{new}}_\mathbb{Q}=(A\otimes
\mathbb{Q},f_*^{\mathrm{new}})$ follows the multiplicative group law $x+y-xy$. By the
universal property of the $K$-theory, there exists a unique morphism of cohomology
theories
$$
\mathrm{ch}\colon K\longrightarrow A^{\mathrm{new}}_\mathbb{Q}
$$
and this is just Grothendieck's Riemann-Roch theorem. In fact this ring morphism
$\mathrm{ch}\colon K(X)\to A (X)\otimes \mathbb{Q}$ is the \textbf{Chern character},
the additive extension of the series $e^t$, because by (\ref{10}) we have
$$
\mathrm{ch}(L_x)=1-c_1^\mathrm{new}(L_x^*)=1-c_1^\mathrm{new}(L_{-x})=1-(1-e^x)=e^x.
$$
If we consider the multiplicative extension $\mathrm{Td}$ of the series
$$
F(t)^{-1}=\frac{t}{1-e^{-t}}=1+\frac{t}{2}+\frac{t^2}{12}-\frac{t^4}{720}+\ldots ,
$$
usually named \textbf{Todd class}, then we obtain\bigskip

\noindent \textbf{Grothendieck's Riemann-Roch Theorem:} \emph{Let $A $ be a
cohomology theory on the smooth quasi-projective varieties over a field following the
additive law. For any projective morphism $f\colon Y\to X$, we have a commutative
square}
$$
\xymatrix{K(Y) \ar[r]^-{f_!} \ar[d]_{\mathrm{Td}(T_Y)\cdot
\mathrm{ch}} & K(X) \ar[d]^{\mathrm{Td}(T_X)\cdot \mathrm{ch}} \\
A (Y)\otimes \mathbb{Q} \ar[r]^-{f_*} &A (X)\otimes \mathbb{Q} }
$$

\noindent \emph{Proof:} Since $\mathrm{ch}\colon K\to A^{\mathrm{new}}_\mathbb{Q}$
preserves direct images,
\begin{align*}
\mathrm{ch}(f_!(y))&=f^\mathrm{new}_*(\mathrm{ch}(y))=f_*\bigl[ F(f^*T_X-T_Y)
\mathrm{ch}(y)\bigr] \\
&=F(T_X)f_*\bigl[ F(T_Y)^{-1}\mathrm{ch}(y)\bigr] =\mathrm{Td}(T_X)^{-1}f_*\bigl[
\mathrm{Td}(T_Y)\mathrm{ch}(y)\bigr] .
\end{align*}
\qed

\medskip

\noindent \textbf{Definition:} A \textbf{graded cohomology theory} is a cohomology
theory with values in the category of graded commutative rings (remark that elements
of negative degree are assumed to be null)
$$
A^\bullet (X)=\mbox{$\bigoplus \limits _{n\ge 0}$}A^n(X),
$$
such that, for any projective morphism $f\colon Y\to X$ between connected smooth
quasi-projective varieties, the direct image $f_*\colon A^n(Y)\to A^{n+d}(X)$ changes
the degree in the codimension $d=\dim X-\dim Y$.

Morphisms of graded cohomology theories are defined to be homogeneous (degree
preserving) morphisms of cohomology theories.\medskip

Remark that the fundamental class of a smooth closed subvariety $Y\to X$ of
codimension $d$ is in $A^d(X)$, that the Chern class of a line bundle $L\to X$ is in
$A^1(X)$ and that, in general, $c_n(E)$ is in $A^n(X)$.\bigskip

\noindent \textbf{Examples:} The graded $K$-theory, the Chow ring and the singular
cohomology ring in the complex case, are graded cohomology theories.

\begin{lemma}
Any graded cohomology theory follows the additive group law
$$
c_1(L_x\otimes L_y)=x+y.
$$
\end{lemma}

\noindent \emph{Proof:} Let us consider the line bundles $\pi _1^*\xi _m$ and $\pi
_2^*\xi _n$ on $\mathbb{P}^m\times \mathbb{P}^n$, and the component $A^0$ of degree 0
in the cohomology ring of a point.

According to axiom 6 we have $A^1(\mathbb{P}^m\times \mathbb{P}^n)=A^0(\pi
_1^*x_m)\oplus A^0(\pi _2^*x_n)$, so that $c_1(\pi _1^*\xi _m\otimes \pi _2^*\xi _n
)=a(\pi _1^*x_m)+b(\pi _2^*x_n)$, and it is easy to show that $a=b=1$.

Now, by Jouanolou's trick we have $p^*L_x=f^*(\pi _1^*\xi _m)$ and $p^*L_y=f^*(\pi
_2^*\xi _n)$ for some morphism $f\colon P\to \mathbb{P}^m\times \mathbb{P}^n$, so
that
$$
p^*c_1(L_x\otimes L_y)=f^*c_1(\pi _1^*\xi _m\otimes \pi _2^*\xi _n)=f^*\bigl( c_1(\pi
_1^*\xi _m)+c_1(\pi _2^*\xi _n)\bigr) =p^*(x+y).
$$

Since $p^*$ is an isomorphism, we conclude.\hfill $\square $

\bigskip

\begin{theorem}\label{universalgk}
Let $A^\bullet $ be a graded cohomology theory on smooth quasi-projective varieties
over a perfect field. There exists a unique morphism of graded cohomology theories
$GK^\bullet \to A^\bullet \otimes \mathbb{Q}$.
\end{theorem}

\noindent \emph{Proof:} Let $X$ be a smooth quasi-projective variety and let $Y\to X$
be a closed subvariety of codimension $d$. If $Y$ is smooth, by the Riemann-Roch
theorem for the closed immersion $Y\to X$, we have
$$
\mathrm{ch}(\mathcal{O}_Y)=\mathrm{ch}(i_!(1))=[Y]+\ldots \in \mbox{$\bigoplus
\limits _{n\ge d}A^n(X)\otimes \mathbb{Q}$}.
$$

In general, since $k$ is perfect, any closed subvariety $Y$ is smooth outside a
closed set $Y_{\mathrm{sing}}$ of bigger codimension.

If we consider the open immersion $j\colon U=X-Y_{\mathrm{sing}}\to X$, then by axiom
3 we have injective morphisms (recall that $A^i(Y_{\mathrm{sing}})=0$ when $i<0$)
$$
j^*\colon A^n(X)\otimes \mathbb{Q}\longrightarrow A^n(U)\otimes \mathbb{Q},\ n\le d.
$$

Since $j^*(\mathrm{ch}(\mathcal{O}_Y))=\mathrm{ch}(j^!\mathcal{O}_Y)=[Y\cap U]+...\in
A^\bullet (U)\otimes \mathbb{Q}$, we see that $\mathrm{ch}\colon K(X)\to A^\bullet
(X)\otimes \mathbb{Q}$ preserves filtrations, $\mathrm{ch}(F^d(X))\subseteq \oplus
_{n\ge d}A^n(U)\otimes \mathbb{Q}$.

Hence it induces a homogeneous ring morphism $\varphi \colon GK^\bullet (X)\to
A^\bullet (X)\otimes \mathbb{Q}$ preserving inverse images and, by Panin's lemma, to
prove that it preserves direct images, we have to show that it preserves the Chern
class of any line sheaf $L$.

Now, if we put $\alpha =c_1(L)\in A^1(X)\otimes \mathbb{Q}$, then
$$
\varphi (c_1^{GK}(L))=\varphi ([L-1])=[\mathrm{ch}(L-1)]=[e^{\alpha }-1]=[\alpha
+\ldots ]=\alpha .
$$

This is the unique possible morphism of graded theories $\varphi \colon GK^\bullet
\to A^\bullet \otimes \mathbb{Q}$, since we have $\varphi ([\mathcal{O}_Y])=\varphi
([Y]^{GK})=[Y]^A$ for any smooth closed subvariety $Y\to X$, and the above argument
shows that this condition fully determines $\varphi ([\mathcal{O}_Y])$ when $Y$ is an
arbitrary closed subvariety of $X$ (recall that these classes $[\mathcal{O}_Y]$
generate $GK^\bullet (X)$ as a group).\qed

\section{Applications}

\begin{enumerate}

\item
Theorem \ref{universalgk} shows that any numerical invariant (global intersection
numbers, Chern numbers,...) defined in $A^\bullet (X)\otimes \mathbb{Q}$ coincides
with the number obtained in $GK^\bullet (X)\otimes \mathbb{Q}$. For example, in the
case of the singular cohomology of a complex projective smooth variety $X$, the
topological self-intersection number of the diagonal in $X\times X$ is just the
topological Euler-Poincar\'{e} characteristic $\chi _{\mathrm{top}}(X)$, so that it
coincides with the algebraic self-intersection number,
$$
\mbox{$\chi _{\mathrm{top}}(X)=\sum \limits _{p,q}(-1)^{p+q}\mathrm{dim}\,
H^p(X,\Omega _X^q)$},
$$
and, in the case of a smooth projective curve $C$ we obtain the coincidence of the
topological genus, $\chi _{\mathrm{top}}(C)=2-2g_{\mathrm{top}}$, with the algebraic
genus $g=\mathrm{dim}\, H^0(C,\Omega _C^1)$.

Moreover, if $X$ is connected of dimension $d$, then $\chi _{\mathrm{top}}(X)$ is the
degree of the topological Chern class $c_d(T_X)$, so that it also coincides with the
degree of the algebraic Chern class $c_d^{GK}(T_X)$, which is just the number of
zeroes of any algebraic tangent vector field with isolated singularities.\bigskip

\item
The Chern character and Todd class of a locally free sheaf $E$ of rank $r$ are
\begin{align*}
\mathrm{ch}(E)&=\mbox{$\sum \limits _ie^{\alpha _i}$}=r+c_1+
\mbox{$\frac{1}{2}(c_1^2-2c_2)+\frac{1}{6}(c_1^3-3c_1c_2+3c_3)$}+\ldots \\
\mathrm{Td}(E)&=\mbox{$\prod \limits _i(1+\frac{\alpha _i}{2}+\frac{\alpha _i^2}{12}+
\ldots )$}=1+\mbox{$\frac{1}{2}$}c_1+\mbox{$\frac{1}{12}$} (c_1^2+c_2)+
\mbox{$\frac{1}{24}$}c_1c_2 +\ldots
\end{align*}
Let $C$ be a smooth projective irreducible curve over an algebraic closed field and
$K\coloneqq c_1(\Omega _C)=-c_1(T_C)\in A^1(C)\otimes \mathbb{Q}$ be the class of a
canonical divisor. Then the Riemann-Roch theorem for the projection $\pi \colon C\to
p$ onto a point gives
$$
\chi (C,E)=\pi _!(E)=\pi _*(\mathrm{Td}(T_C)\mathrm{ch}(E))=\mathrm{deg}\,
c_1(E)-\mbox{$\frac{r}{2}$}\mathrm{deg}\, K.
$$
If apply it to $E=\mathcal{O}_C$ we obtain
$$
1-g=\chi (C,\mathcal{O}_C)= -\frac{1}{2}\mathrm{deg}\, K.
$$
If we take $E=L_{D}$ then $c_1(L_D)=[D]$ so that $\mathrm{deg}\, c_1(L_D) =
\mathrm{deg} \, D$. We obtain
$$
\dim H^0(C,L_D)-\dim H^1(C,L_D)= \mathrm{deg}\,D  +1 - g.
$$
By Serre's duality $\dim H^1(C,L_D)=\dim H^0(C,L_{K-D})$ and we obtain the classic
formula for curves.

If $S$ is a smooth projective surface, $K\coloneqq c_1(\Omega ^2_S)=-c_1(T_S)\in
A^1(S)\otimes \mathbb{Q}$ is the class of a canonical divisor, and we put $\chi
_{\mathrm{top}}\coloneqq \mathrm{deg}\, c_2(T_S)$. The Riemann-Roch theorem for the
projection $\pi \colon C\to p$ onto a point gives
\begin{align*}
\chi (S,E)&=\pi _!(E)=\pi _*(\mathrm{Td}(T_S)\mathrm{ch}(E))\\
&=\mbox{$\frac{1}{12}$}r(K^2+\chi _{\mathrm{top}})-\mbox{$\frac{1}{2}$}K\cdot
c_1(E)+\mbox{$\frac{1}{2}$}c_1(E)^2-\mathrm{deg}\, c_2(E).
\end{align*}

When $E$ is the trivial line bundle, we obtain \textbf{Noether's equality}
$$
\chi (S,\mathcal{O}_S)=\mbox{$\frac{1}{12}$}(K^2+\chi _{\mathrm{top}}).
$$

\item
In the case of a closed smooth hypersurface $i\colon Y\to X$, the Riemann-Roch
theorem gives
\begin{align*}
i_*(\mathrm{Td}(T_Y))&=\mathrm{ch}(\mathcal{O}_Y)\cdot \mathrm{Td}(T_X)=(1-e^{-Y})
\cdot \mathrm{Td}(T_X)\\
Y-\mbox{$\frac{1}{2}$}i_*(K_Y)+\ldots &=(Y-\mbox{$\frac{1}{2}$}Y^2+\ldots )
(1-\mbox{$\frac{1}{2}$}K_X+\ldots )
\end{align*}
and we obtain the adjunction formula $i_*(K_Y)=Y(K_X+Y)$.

\item
In the case of a closed smooth subvariety $i\colon Y\to X$ of codimension $d$, the
Riemann-Roch theorem gives
$\mathrm{ch}(\mathcal{O}_Y)=i_*(\mathrm{Td}(N_{Y/X}))=Y+\ldots $. Since
$$
\alpha _1^r+\ldots +\alpha _n^r=p_r(s_1,\ldots ,s_{r-1})+(-1)^{r+1}rs_r
$$
for some polynomial $p_r$, where $s_i$ is the $i$-th elementary symmetric function of
$\alpha _1,\ldots ,\alpha _n$ (just consider the roots of the polynomial  $x^r-1$),
we see that
$$
c_i(\mathcal{O}_Y)=0 \  \forall \, i <d \ , \ c_d(\mathcal{O}_Y)=(-1)^d(d-1)!\cdot Y
$$
modulo torsion; that is to say, in $A^\bullet (X)\otimes \mathbb{Q}$.

\item
If $\omega $ is a rational 1-form on a connected smooth projective surface $S$, let
us determine $c_2(T_S)\in A^2(S)\otimes \mathbb{Q}$ in terms of the singularities of
$\omega $. Let $\Sigma $ be the field of rational functions on $S$. Then
$$
(\omega )(U)=\{ f\omega \in \Omega ^1_S(U)\colon f\in \Sigma \}
$$
is a line sheaf, and $(\omega )\simeq L_D$, where $D$ is the divisor of zeros and
poles of the differential form $\omega $. We have an exact sequence
$$
0\longrightarrow (\omega)\longrightarrow \Omega ^1_S\xrightarrow{\ \wedge \omega \
}\Omega ^2_S\otimes L_{-D}\longrightarrow \mathfrak{C}\longrightarrow 0\leqno{(*)}
$$
If $x,y$ are parameters at a point $p$, and $\omega =h(f\d x+g\d y)$, where $f,g\in
\mathcal{O}_p$ have no common factor, a local equation of $D$ is $h=0$, and
$\mathfrak{C}_p=\mathcal{O}_p/(f,g)$. Now, the Riemann-Roch theorem for the immersion
$p\hookrightarrow S$ gives that $\mathrm{ch}(\mathcal{O}_S/\mathfrak{m}_p)=p$, where
$\mathfrak{m}_p$ denotes the maximal ideal of $p$; hence we have
$$
\mathrm{ch}(\mathfrak{C})=\mbox{$\sum _pl(\mathfrak{C}_p)\cdot p$},
$$
where $l(\mathfrak{C}_p)$ stands for the length of the $\mathcal{O}_{S,p}$-module
$\mathfrak{C}_p$.

If $K\coloneqq c_1(\Omega ^2_S)=c_1(\Omega ^1_S)\in A^1(S)\otimes \mathbb{Q}$ denotes
the class of a canonical divisor, the exact sequence $(*)$ gives the
\textbf{Zeuthen-Segre invariant}:
\begin{align*}
&\mathrm{ch}(L_D)+\mathrm{ch}(L_{K-D})=\mathrm{ch}(\Omega ^1_S)+\mathrm{ch}
(\mathfrak{C}),\\
&e^D+e^{K-D}=2+K+\mbox{$\frac{1}{2}$}K^2-c_2(T_S)+\mathrm{ch}(\mathfrak{C}),\\
&c_2(T_S)=D(K-D)+\mbox{$\sum _pl(\mathfrak{C}_p)\, p$}.
\end{align*}

\end{enumerate}


\begin{thebibliography}{8.}

\bibitem{Borel}
A. Borel and J.P. Serre: {\it Le th\'{e}or\`{e}me de Riemann-Roch}; Bull. Soc. Math. France,
86, 97-136 (1958)

\bibitem{Deglise}
F. D\'{e}glise: {\it Orientation theory in arithmetic geometry}; arXiv:1111.4203

\bibitem{Deglise2}
F. D\'{e}glise: {\it Around the Gysin triangle II}; Doc. Math. 13, 613–675 (2008)

\bibitem{Fulton}
W. Fulton: {\it Intersection theory}; Springer-Verlag, Heidelberg (1998)

\bibitem{SGA6}
A. Grothendieck et al.: {\it Th\'{e}orie des Intersections et Th\'{e}or\`{e}me de Riemann-Roch
(SGA 6)}; Lecture Notes in Math. {\bf 225}, Springer-Verlag, Heidelberg (1971)

\bibitem{Groth}
A. Grothendieck: {\it La th\'{e}orie des classes de Chern}; Bull. Soc. Math. France, 86,
137–154 (1958)

\bibitem{Jouanolou}
J. Jouanolou: {\it Une suite exacte de Mayer-Vietoris en $K$-th\'{e}orie alg\'{e}brique};
Algebraic $K$-theory I (Seattle 1972), Lecture Notes in Math. {\bf 341},
Springer-Verlag, Heidelberg (1973)

\bibitem{Levine}
M. Levine, F. Morel: {\it Algebraic Cobordism}; Springer Monographs in Math.,
Springer-Verlag, Heidelberg (2007)

\bibitem{Milnor}
J. Milnor, J. Stasheff: {\it Characteristic Classes}; Annals Math. Studies {\bf 76},
Princeton Univ. Press, Princeton (1974)

\bibitem{Navarro}
J. A. Navarro: {\it Notes for a Degree in Mathematics};

\url{http://matematicas.unex.es/\textasciitilde navarro/degree.pdf}

\bibitem{Panin}
I. Panin: {\it Riemann-Roch theorems for oriented cohomologies}; Axiomatic, enriched
and motivc homotopy theory, NATO Sci. Ser. II Math. Phys. Chem., vol. 131, Kluwer
Acad. Publ., Dordecht, 261-333 (2004)

\bibitem{Panin2}
I. Panin, K. Pimenov and O. R\"{o}ndigs: {\it A universality theorem for Voevodsky's
algebraic cobordism spectrum}; Homology, Homotopy Appl. 10, no. 2, 211-226 (2008)

\bibitem{Panin3}
I. Panin, K. Pimenov and O. R\"{o}ndigs: {\it On the relation of Voevodsky's algebraic
cobordism to Quillen's $K$-theory}; Invent. Math. 175, 435-451 (2009)

\bibitem{Quillen}
D. G. Quillen: {\it Elementary proofs of some results of cobordism theory using
Steenrod operations}; Advances in Math. 7, 29-56 (1971)

\bibitem{Serre}
J. P. Serre: {\it Alg\`{e}bre Locale. Multiplicit\'{e}s}; Lecture Notes in Math. {\bf 11},
Springer-Verlag, Heidelberg (1965)

\bibitem{Weibel}
C. Weibel: {\it The $K$-book: an introduction to algebraic $K$-theory}; Graduate
Studies in Math. vol. 145, Amer. Math. Soc. (2013)

\end{thebibliography}
\end{document}